\documentclass[a4paper,english,12pt,reqno]{amsart}
\usepackage{babel}
\usepackage{a4wide}
\usepackage{amsmath}
\usepackage{amsthm,amssymb}
\usepackage[all]{xy}
\usepackage{enumerate}
\CompileMatrices

\theoremstyle{plain}
\newtheorem{theorem}{Theorem}[section]
\newtheorem{proposition}[theorem]{Proposition}
\newtheorem{corollary}[theorem]{Corollary}
\newtheorem{lemma}[theorem]{Lemma}

\theoremstyle{definition}
\newtheorem{example}[theorem]{Example}
\newtheorem{definition}[theorem]{Definition}

\newtheorem{remark}[theorem]{Remark}

\newcommand{\BbbZ}{{\mathbb{Z}}}

\DeclareMathSymbol{\Bbbk}{\mathalpha}{AMSb}{'174}

\DeclareMathOperator{\im}{im} 
 
\DeclareMathOperator{\Rat}{Rat} 
\DeclareMathOperator{\Hom}{Hom}

\newcommand{\tensor}[1][]{\otimes_{#1}}

\newcommand{\lMod}[1]{{}_{#1}\mathcal{M}}
\newcommand{\rMod}[1]{\mathcal{M}_{#1}}
\newcommand{\ModHom}[1]{\Hom_{#1}}

\newcommand{\rCom}[1]{\mathcal{M}^{#1}}
\newcommand{\ComHom}[1]{\Hom^{#1}}

\newcommand{\Pleft}{\ensuremath{\mathbf{P_\ell}}}
\newcommand{\PellAlg}[1][R]{\mathrm{P_\ell Alg}_{#1}}

\newcommand{\bil}[2]{\langle #1, #2 \rangle}

\begin{document}

\title[Hopf Algebras over Commutative Rings]{Duality and Rational Modules in Hopf Algebras over Commutative Rings.}
\author{J. Y. Abuhlail}
\address{Institute of Mathematics. University of D\"{u}sseldorf. D-40225 D\"{u}sseldorf, Germany}
\email{abuhlail@math.uni-duesseldorf.de}
\author{J. G\'{o}mez-Torrecillas}
\address{Departamento de \'{A}lgebra. Universidad de Granada. E-18071 Granada, Spain}
\email{torrecil@ugr.es}
\author{F. J. Lobillo}
\address{Departamento de \'{A}lgebra. Universidad de Granada. E-18071 Granada, Spain}
\email{jlobillo@ugr.es}
\thanks{The authors are grateful for the financial support by the
  Spanish--German exchange program \emph{Acciones Integradas}.}
\date{\today}
\begin{abstract}
Let $A$ be an algebra over a commutative ring $R$. If $R$ is noetherian and
$A^\circ$ is pure in $R^A$, then the categories of rational left $A$--modules
and right $A^\circ$--comodules are isomorphic. In the Hopf algebra case,
we can also strengthen the Blattner-Montgomery duality theorem. Finally, we
give sufficient conditions to get the purity of $A^\circ$ is $R^A$.
\end{abstract}

\maketitle

\section*{Introduction}

It is well known that the theory of Hopf algebras over a field cannot be
trivially passed to Hopf algebras over a commutative ring. For instance let us
consider $\BbbZ[x]$ as Hopf algebra and let $\mathfrak{a}$ be the Hopf 
ideal generated by $\langle 4,2x \rangle$. Let $H$ be the Hopf
$\BbbZ$--algebra $H = \BbbZ[x] / \mathfrak{a}$. The finite dual is zero in
this situation. However $H \cong \BbbZ_4[x] / \langle 2x \rangle$, so we can
view $H$ as a Hopf $\BbbZ_4$--algebra. If $I$ is a $\BbbZ_4$--cofinite ideal
of $H$ then every element nonzero in $H/I$ has order $2$, which implies every
element in $H^\circ$ has order $2$. In this situation $H^\circ$ is not pure in
$\BbbZ_4^H = \operatorname{Map}(H,\BbbZ_4)$ as a $\BbbZ_4$--module (since
$H^\circ$ is not free) and a canonical $\BbbZ_4$--coalgebra structure on
$H^\circ$ cannot be expected. 

A basic result from the theory of coalgebras over fields is that the comodules
are essentially rational modules. Thus, given a (left) non-singular pairing of
a coalgebra $C$ and an algebra $A$ (see \cite{Radford:1973}), the categories
of rational left $A$--modules and right $C$--comodules are isomorphic. This
applies in particular for the canonical pairings $(C,C^*)$ and $(A^{\circ},A)$
derived from a coalgebra $C$ and an algebra $A$, respectively. An attempt to
develop systematically the theory of rational modules associated to a pairing
$(C,A)$, where $C$ is a coalgebra and $A$ is an algebra over an arbitrary
commutative ring $R$, is \cite{Gomez:1998}. A corollary of the theory there
developed is that if $C$ is projective as an $R$--module, then the category of
right $C$--comodules is isomorphic to the category of rational left modules
over the convolution dual $R$--algebra $C^{*} = \ModHom{R}(C,R)$ (this result
was independently obtained in \cite{Wisbauer:1997} by means of a different
approach). However, the results of \cite{Gomez:1998} are proved in a framework
which does not allow to apply them directly to the pairing $(A^{\circ},A)$, for
a given $R$--algebra $A$. In fact, the first problem is to endow the finite
dual $A^{\circ}$ with a comultiplication, which entails some serious technical
difficulties at the very beginning due to the lack of exactness of the tensor
product bifunctor $-\otimes_R-$. Nevertheless, it has been recently proved in
\cite[Theorem 2.8]{Abuhlail/Gomez/Wisbauer:2000} that if $R$ is noetherian and
$A^{\circ}$ is pure in the $R$--module $R^A$ of all maps from $A$ to $R$, then
$A^{\circ}$ is a coalgebra. We have observed that the notion of rational
pairing introduced in \cite{Gomez:1998} can be restated in order that the
methods developed there can be applied to the pairing $(A^{\circ},A)$ to prove
that the category of right $A^{\circ}$--comodules is isomorphic to the
category of rational left $A$--modules. This applies, in particular, for any
algebra $A$ over a hereditary noetherian commutative ring. We explain our
general theory of rational modules and comodules in Sections 2 and 3.

We apply our methods to strengthen the Blattner-Montgomery duality
theorem for Hopf algebras over commutative rings. Let $H$ be a Hopf
algebra over the commutative ring $R$. When $R$ is a field, the
Blattner--Montgomery duality theorem says that if $U$ is a Hopf
subalgebra of $H^\circ$, $A$ is an $H$--module algebra such that
the $H$--action is locally finite in a sense appropriate to the
choice of $U$, $H$ and $U$ have bijective antipodes and there is a
certain right-left symmetry in the action of $H \# U$ on $H$, then
\begin{displaymath}
(A \# H) \# U \cong A \tensor (H \# U)
\end{displaymath}

There are two proofs of this theorem in the literature. The first one appeared
in \cite[Theorem 2.1]{Blattner/Montgomery:1985}, and a new one, due to Blattner,
appears in \cite[Theorem 9.4.9]{Montgomery:1993}. Since, in this situation, the
$U$--comodules are just the $U$--locally finite $H$--modules, it is easy to see
that the two theorems are equivalent. In the case of a general commutative ring
$R$, there is a similar theorem due to Van den Bergh (see
\cite{VandenBergh:1984}) when $H$ is finitely generated and projective over $R$.
A generalization of \cite{Blattner/Montgomery:1985} for Hopf algebras over a
Dedekind domain $R$ was proved by Chen and Nichols, under the technical
condition that $U$ is $R$--closed in $H^\circ$ (see \cite[Theorem
5]{Cao-Yu/Nichols:1990}). This condition guarantees that every $U$--locally
finite is rational. However, it is not evident that \cite[Theorem
5]{Cao-Yu/Nichols:1990} generalizes \cite[Theorem 9.4.9]{Montgomery:1993}.

We show that the ideas used in the proof of \cite[Theorem
9.4.9]{Montgomery:1993}, together with our results on rational modules and
comodules, can be combined to get a duality theorem for Hopf algebras over a
noetherian commutative ring $R$ which generalizes both \cite[Theorem
9.49]{Montgomery:1993} and \cite[Theorem 5]{Cao-Yu/Nichols:1990} and, hence,
\cite[Theorem 2.1]{Blattner/Montgomery:1985}.

In Section 4 we introduce a class of $R$-algebras $\PellAlg$ (in
case $R$ is noetherian) which satisfy the property that $A^{\circ }
\subset R^{A}$ is an $R$-pure submodule for each $A \in \PellAlg$
and, hence, the canonical pairing $(A^{\circ },A, \bil{-}{-})$ is a
rational pairing. We give several examples of such algebras, among
them the polynomial algebra $R[x_{1},\dots,x_{n}]$ and the algebra
of Laurent polynomials
$R[x_{1},x_{1}^{-1},\dots,x_{n},x_{n}^{-1}]$.

\section{Preliminaries and Basic Notions}

In this paper $R$ is a commutative ring with unit. Let $A$ be an
associative $R$--algebra with unit. The category of all left
$A$--modules is denoted by $\lMod{A}$. As usual, the notation $X
\in \mathcal{C}$, where $\mathcal{C}$ is a category, means $X$ is
an object of $\mathcal{C}$. By $\tensor$ we denote the tensor
product $\tensor[R]$ unless otherwise explicitly stated. Moreover,
if $\pi \in S_n$ (the symmetric group on $n$ symbols) then
$\tau_\pi$ is the canonical isomorphism
\begin{equation*}
\tau_\pi : M_1 \tensor \dots \tensor M_n \longrightarrow M_{\pi(1)} \tensor
\dots \tensor M_{\pi(n)}
\end{equation*}
Let $M,X$ be $R$--modules. If $N$ is an $R$--submodule of $M$ then $N$ is
called $X$--pure if $N \tensor X \subseteq M \tensor X$. The inclusion $N 
\subseteq M$ is called pure if $N$ is $X$--pure for all $R$--modules
$X$. Unless otherwise stated, pure, projective and flat mean pure, projective
and flat in $\lMod{R}$.

Let $\mathcal{A}$ be a Grothendieck category. A \emph{preradical} for
$\mathcal{A}$ is a subfunctor of the identity endofunctor $id_{\mathcal{A}} :
\mathcal{A} \rightarrow \mathcal{A}$. We follow \cite{Stenstrom:1975} for
categorical basic notions.

\begin{definition}
An $R$--coalgebra is an $R$--module $C$ together with two homomorphisms of
$R$--modules
\begin{equation*}
\Delta: C \rightarrow C \tensor C \text{ (comultiplication) and } \epsilon : C
\rightarrow R \text{ (counit)}
\end{equation*}
such that
\begin{gather*}
(id_C \tensor \Delta) \circ \Delta = (\Delta \tensor id_C) \circ \Delta \quad
\text{and} \\ (id_C \tensor \epsilon) \circ \Delta = (\epsilon \tensor id_C)
\circ \Delta = id_C
\end{gather*}
\end{definition}

\begin{definition}
A right $C$--comodule is an $R$--module together with an $R$--homomorphism
\begin{equation*}
\rho_M : M \rightarrow M \tensor C
\end{equation*}
such that
\begin{gather*}
(id_M \tensor \Delta) \circ \rho_M = (\rho_M \tensor id_C) \circ \rho_M \quad
\text{and} \\ (id_M \tensor \epsilon) \circ \rho_M = id_M.
\end{gather*}

Let $M,N$ be right $C$--comodules. A homomorphism of $R$--modules $f : M
\rightarrow N$ is said to be a comodule morphism (or $C$--colinear) if $\rho_N
\circ f = (f \tensor id) \circ \rho_M$. By $\ComHom{C}(M,N)$ we denote the
$R$--module of all colinear maps from $M$ to $N$. The right $C$--comodules with
the $C$--colinear maps between them constitute an additive category denoted by
$\rCom{C}$. In case that $C$ is a flat $R$--module, $\rCom{C}$ is a
Grothendieck category (see \cite[Corollary 3.15]{Wisbauer:1997})
\end{definition}

For basic notions on coalgebras and comodules over commutative rings we refer to
\cite{Gomez:1998,Wisbauer:1997}.

\section{Rational modules}

In \cite{Gomez:1998}, the theory of rational modules is developed under the
assumption of the existence of the there--called rational pairing. The main
example exhibited there was $(C^*,C)$ with $C$ an $R$--projective coalgebra. Our
aim is to deal with the finite dual of an $R$--algebra, so in this section we
provide a weaker definition of rational pairing (in order to cover the finite
dual example) which also implies the results on rational modules developed in
\cite{Gomez:1998}.

\subsection{Rational Systems}

Let $A,P$ be $R$--modules and let
\begin{displaymath}
\bil{-}{-}: P \times A \rightarrow R
\end{displaymath}
be a bilinear form. For every $R$--module $M$, define the $R$--linear map
\begin{equation*}
\begin{split}
\alpha_M : M \tensor P &\longrightarrow \ModHom{R}(A,M) \\ m \tensor p
&\longmapsto \big[ a \mapsto m \bil{p}{a} \big]
\end{split}
\end{equation*}

\begin{proposition}\label{prop-q}
In the previous situation the following statements are equivalent:
\begin{enumerate}
\item $\alpha_M$ is injective
\item If $\sum m_{i} \tensor p_{i} \in M \tensor P$, then $\sum m_{i} \tensor p_{i} = 0$ if and only if for every $a \in A$, $\sum m_{i} \bil{p_{i}}{a} = 0$.
\end{enumerate}
\end{proposition}

\begin{proof}
Note that for each $R$--module $M$,
\begin{equation*}
\ker(\alpha_{M}) = \{\sum m_{i} \tensor p_{i}~|~\sum \bil{p_{i}}{-}m_{i}=0\}.
\end{equation*}
\end{proof}

\begin{definition}\label{proQ}
The three-tuple $(P,A,\bil{-}{-})$ is a \emph{rational system} if
$\alpha_M$ is injective for every $R$--module $M$.
\end{definition}

\begin{remark}
By \cite[Proposition 2.3]{Gomez:1998} and Proposition \ref{prop-q}, a rational
system as defined in \cite[Definition 2.1]{Gomez:1998} is a rational system in
the present setting.
\end{remark}

\begin{remark}\label{flat}
Let $(P,A,\bil{-}{-})$ be a rational system. Let $M$ be an $R$--module and let
$N$ be an $R$--submodule of $M$. Consider the following commutative diagram
\begin{equation*}
\xymatrix{ N \tensor P \ar@{^{(}->}[r]^-{\alpha_N} \ar[d]_{i_N \tensor id_P} &
\ModHom{R} (A,N) \ar@{^{(}->}[d]^{i} \\ M \tensor P \ar@{^{(}->}[r]_-{\alpha_M}
& \ModHom{R} (A,M) }
\end{equation*}
Note that $\alpha_{N}$ is injective since the
three-tuple $(P,A,\bil{-}{-})$ is a rational system, and
$i:\ModHom{R}(A,N) \rightarrow \ModHom{R}(A,M)$, $g \mapsto i_{N}
\circ g$ is injective. Hence $i_{N} \tensor id_{P}$ is injective.
Since $M$ was arbitrary in $\lMod{R}$ we conclude that $P$ should
be flat as an $R$--module.
\end{remark}

The following Proposition replaces \cite[Proposition 2.2]{Gomez:1998} in order 
to show that the canonical comodule structure over a rational module is
pseudocoassociative 

\begin{proposition}\label{q-2}
If $(P,A,\bil{-}{-})$ and $(Q,B,\bil{-}{-})$ are rational systems, then the
induced pairing
\begin{equation*}
[-,-]: P \tensor Q \times A \tensor B \rightarrow R
\end{equation*}
defined by
\begin{equation*}
[ \sum p_{i} \tensor q_{i},\sum a_{j} \tensor b_{j} ] = \sum \bil{p_{i}}{a_{j}}
\bil{q_{i}}{b_{j}}
\end{equation*}
is a rational system.
\end{proposition}

\begin{proof}
Consider the commutative diagram
\begin{equation*}
\xymatrix{ M \tensor P \tensor Q \ar[r]^{\beta_M} \ar[d]^{\alpha_{M \tensor P}}
& \ModHom{R}(A \tensor B, M) \ar[d]^{\eta} \\ \ModHom{R} (B, M \tensor P)
\ar[r]^-{(\alpha_M)_*} & \ModHom{R} (B, \ModHom{R}(A,M))}
\end{equation*}
where $\eta$ is the adjunction isomorphism, $\beta_{M}$ is the mapping analogous
to $\alpha$ with respect to the pairing $[-,-]$ and $(\alpha_{M})_*$ is the
homomorphism induced by $\alpha_{M}$. This last morphism is monic and
$\alpha_{M\tensor P}$ is monic. Therefore, $\beta_{M}$ is a monomorphism.
\end{proof}

For each set $S$, let $R^S$ denote the $R$--module of all maps from $S$ to $R$.
If $\alpha_R$ is injective, then $P$ is isomorphic to a submodule of $A^* =
\ModHom{R} (A, R) \subseteq R^A$. We identify the $R$--module $P$ with its image
in $R^A$, so every $p \in P$ is identified with the $R$--linear map $\bil{p}{-}
: A \rightarrow R$.

\begin{definition}\label{mockproj}
We say $P$ is \emph{mock-projective} (relative to $A$ and
$\bil{-}{-}$) if $\alpha_R$ is injective and for every $p_1, \dots,
p_n \in P$ there are $a_1, \dots, a_m \in A$ and $g_1, \dots, g_m
\in R^A$ such that for every $i = 1, \dots, n$, $p_i = \sum \bil{p_i}{a_l}
g_l$.
\end{definition}

Proposition \ref{prop-q} can be improved under the assumption that $P$ is mock-projective:

\begin{proposition}\label{prop-Q}
Assume $P$ is mock-projective. If $M \in \lMod{R}$, then $P
\subseteq R^A$ is $M$--pure if and only if $\alpha_M$ is injective.
Therefore, $(P,A,\bil{-}{-})$ is a rational system if and only if
$P \subseteq R^A$ is pure.
\end{proposition}

\begin{proof}
Assume $P \subseteq R^A$ is $M$--pure and let $\sum m_i \tensor p_i \in
M \tensor P$. Assume $\sum \bil{p_i}{a}m_i = 0$ for all $a \in A$. By
Definition \ref{mockproj} there are $a_1, \dots, a_m \in A$ and $g_1, \dots,
g_m \in R^A$ such that
\begin{equation*}
p_i = \sum \bil{p_i}{a_l} g_l, \text{ for each $i=1, \dots, n$}.
\end{equation*}
As $\sum m_i \tensor p_i \in M \tensor R^A$  we have
\begin{equation}\label{eq:escero}
\begin{split}
\sum m_i \tensor p_i &= \sum m_i \tensor \sum \bil{p_i}{a_l} g_l = \sum m_i \bil{p_i}{a_l} \tensor g_l \\
&= \sum \big( \sum m_i \bil{p_i}{a_l} \big) \tensor g_l = \sum 0 \tensor g_l = 0,
\end{split}
\end{equation}
and by purity $\sum m_i \tensor p_i = 0$ as an element in $M \tensor P$. By
Proposition \ref{prop-q} $\alpha_M$ is injective.

Assume now $\alpha_M$ is injective. The diagram
\begin{equation*}
\xymatrix{ M \tensor P \ar[d]_-{id \tensor i} \ar[r]^-{\alpha_M} &
\ModHom{R}(A,M) \ar[d]^-{i} \\ M \tensor R^A \ar[r]_-{\beta_M} & M^A }
\end{equation*}
is commutative where $i$ denotes inclusion and $\beta_M(m \tensor f)(a) = m
f(a)$. Since $i \circ \alpha_M$ is injective, we get $id \tensor i$ is injective
and $P \subseteq R^A$ is $M$--pure.
\end{proof}

\subsection{Rational Pairings}
Let $(C,\Delta,\epsilon)$ be an $R$--coalgebra and let $(A,m,u)$ be
an $R$--algebra.

\begin{definition}
A \emph{rational pairing} is a rational system $(C,A,\bil{-}{-})$
where $C$ is an $R$--coalgebra, $A$ is an $R$--algebra and the map
$\varphi : A \rightarrow C^*$ given by $\varphi(a)(c) = \bil{c}{a}$
is a homomorphism of $R$--algebras. This is equivalent to require
\begin{equation}\label{eq:bilphi}
\begin{gathered}
\bil{-}{-} \circ (id_C \tensor m) = (\bil{-}{-} \tensor \bil{-}{-})
\circ \tau_{(23)} \circ (\Delta \tensor id_A^{\tensor 2}) \\
\text{and} \quad \epsilon = \bil{-}{-} \circ (id_C \tensor u) =
\bil{-}{1}.
\end{gathered}
\end{equation}

\end{definition}

Now we can parallel the definition and properties of rational modules in
\cite{Gomez:1998}. The proofs are formally the same as in \cite{Gomez:1998}.
We include some of them for convenience of the reader.

\begin{definition}
Let $T = (C,A,\bil{-}{-})$ be a rational pairing. An element $m$ in a left
$A$--module $M$ is called \emph{rational} (with respect to the pairing $T$) if
there exist finite subsets $\{m_i\} \subseteq M$ and $\{c_i\} \subseteq C$ such
that $a m = \sum m_i \bil{c_i}{a}$ for every $a \in A$. The subset
$\{(m_i,c_i)\} \subseteq M \times C$ is called a \emph{rational set of
parameters for $m$} (with respect to the pairing $T$). The subset $\Rat^T(M)$ of
$M$ consisting of all rational elements of $M$ is clearly an $R$--submodule of
$M$. A left $A$--module is called \emph{rational} (with respect to the pairing
$T$) if $M = \Rat^T(M)$. The full subcategory of $\lMod{A}$ whose objects are
all the rational (with respect to the pairing $T$) left $A$--modules will be
denoted by $\Rat^T(\lMod{A})$. We use the notation $\Rat$ instead of $\Rat^T$
when the rational pairing $T$ is clear from the context.
\end{definition}

\begin{remark}\label{rattensor}
As a consequence of Proposition \ref{prop-q}, if $\{(m_i,c_i)\} \subseteq M
\times C$ is a rational set of parameters for $m$ and $\sum m_i \tensor c_i =
\sum n_j \tensor d_j$, then $\{(n_j,d_j)\} \subseteq M \times C$ is a rational
set of parameters for $m$. In fact a rational set of parameters for $m$ can be
viewed as a representative of an element $\sum m_i \tensor c_i \in M \tensor
C$.
\end{remark}

Following \cite{Wisbauer:1991}, we will denote by $\mathbf{\sigma}[{}_{A}C]$
the full subcategory of $\lMod{A}$ consisting of all the left $A$--modules
subgenerated by ${}_{A}C$. This means that a left $A$--module belongs to
$\mathbf{\sigma}[{}_{A}C]$ if and only if it is isomorphic to a submodule of a
factor module of a direct sum of copies of ${}_{A}C$.

\begin{theorem}\label{iso}%\label{sigma}
Let $T = (C,A,\bil{-}{-})$ be a rational pairing. Then
\begin{enumerate}
\item $\Rat = \Rat^T : \lMod{A} \rightarrow \lMod{A}$ is a left exact
  preradical.
\item The categories $\Rat(\lMod{A})$ and $\rCom{C}$ are isomorphic.
\item $\Rat(\lMod{A}) = \mathbf{\sigma}[{}_{A}C]$.
\end{enumerate}
\end{theorem}

\begin{proof}
The proofs of these facts are formally the same that \cite[Propositions 2.9
and 2.10, Theorems 3.12 and 3.13]{Gomez:1998}, using Propositions \ref{prop-q} 
and \ref{q-2} instead of \cite[Propositions 2.3 and 2.2]{Gomez:1998}
\end{proof}

The isomorphism given in Theorem \ref{iso} is defined in terms of sets of rational parameters: if $M \in \Rat(\lMod{A})$ then the structure of right $C$
comodule is $\omega_M(m) = \sum m_i \tensor c_i$ where
$\{(m_i,c_i)\}$ is a set of rational parameters, if $(M,\delta_M) \in
\rCom{C}$ and $\delta_M(m) = \sum m_i \tensor c_i$, then $\{(m_i,c_i)$ is a
set of rational parameters for $m$. See \cite[Propositions 3.5 and
3.11]{Gomez:1998} for details.

Sweedler's $\Sigma$--notation can be introduced in terms of sets of rational
parameters. Let $T = (C,A,\bil{-}{-})$ be a rational pairing. We have the
injective $R$--linear map $\alpha_R : C \rightarrow A^*$ defined by
\begin{equation*}
\begin{split}
\alpha_R : C &\longrightarrow A^* \\ c &\longmapsto \big[ a \mapsto \bil{c}{a}
\big]
\end{split}
\end{equation*}
Let us regard $A^*$ as a left $A$--module via $(a \lambda)(b) = \lambda(ba)$
and let us identify $C$ with $\alpha_R(C)$. As in \cite[Proposition
3.2]{Gomez:1998}, $C = \Rat({}_A A^*)$. Note that the set of rational
parameters for $c \in C$ is given by $\Delta(c)$. If $\{(c_1,c_2)\}_{(c)} =
\{(c_1,c_2)\} \subseteq C \times C$ represents a set of rational parameters of
$c \in C$ then the comultiplication can be represented as 
\begin{equation*}
\Delta(c) = \sum_{(c)} c_1 \tensor c_2 = \sum c_1 \tensor c_2
\end{equation*}
Note that \eqref{eq:bilphi} means
\begin{equation*}
\bil{c}{ab} = \sum \bil{c_1}{a} \bil{c_2}{b}
\end{equation*}
in $\Sigma$--notation. Analogously, Let $(M,\delta_M) \in \rCom{C}$. The set
of rational parameters for $m \in M$ is given by $\delta_M(m)$. We are going
to use Sweedler's $\Sigma$--notation on $C$--comodules, i.e., $\delta_M(m) =
\sum_{(m)} m_0 \tensor m_1$ where $\{(m_0,m_1)\}_{(m)} \subseteq M \times C$
represents an arbitrary set of rational parameters for $m \in M$.

\subsection{The finite dual coalgebra $A^\circ$}

\emph{In this subsection, the commutative ring $R$ is assumed to be
noetherian}. Let $A$ be an $R$--algebra. Recall that the canonical
structure of an $A$--bimodule on $A^*$ is given by
\begin{equation} \label{eq:dual-co}
(af)(b)=f(ba)\text{ and }(fa)(b)=f(ab)\text{ for }f\in A^*\text{ and }a,b\in A.
\end{equation}
Let $A^\circ = \{ f \in A^*~|~\text{$Af$ is finitely generated as an
$R$--module}\}$. Then by \cite[Proposition 2.6]{Abuhlail/Gomez/Wisbauer:2000},
\begin{equation}\label{eq:dual-co2}
\begin{split}
A^\circ &= \{ f \in A^*~|~\text{$Af$ is finitely generated as an $R$--module}\}
\\ &= \{ f \in A^*~|~\text{$fA$ is finitely generated as an $R$--module}\} \\ &=
\{ f \in A^*~|~\text{$\ker f$ contains an $R$--cofinite ideal of $A$}\} \\ &= \{
f \in A^*~|~\text{$\ker f$ contains an $R$--cofinite left ideal of $A$}\} \\ &=
\{ f \in A^*~|~\text{$\ker f$ contains an $R$--cofinite right ideal of $A$}\}.
\end{split}
\end{equation}
By \cite[2.3]{Abuhlail/Gomez/Wisbauer:2000}, $A^\circ$ is an $A$--subbimodule of
$A^*$.

\begin{remark}
If $A$ is finitely generated and projective in $\lMod{R}$ then $A^\circ = A^*$
is pure in $R^{A}$. In this case $(A \tensor A)^* \simeq A^* \tensor A^*$.
\end{remark}

By \cite[Theorem 2.8]{Abuhlail/Gomez/Wisbauer:2000}, $A^\circ$ is an
$R$--coalgebra whenever $A^{\circ}$ is a pure submodule of $R^A$.
In this section we will prove that $(A^\circ,A,\bil{-}{-})$ is a
rational pairing. So we are going to describe the right
$A^\circ$--comodules as rational left $A$--modules. This applies in
particular when $R$ is hereditary.

Next lemma is used to prove the rationality of the three-tuple
$(A^\circ,A,\bil{-}{-})$.

\begin{lemma}\label{movemaps}
Let $S$ be a set and let $f_1, \dots, f_n \in R^S$. Then there
exist $s_1, \dots, s_m \in S$ and $g_1, \dots, g_m \in R^S$ such
that
\begin{equation*}
f_i = \sum f_i(s_l) g_l, \text{ for each $i=1, \dots, n$}.
\end{equation*}
In particular, if $A,P$ are $R$--modules and $\bil{-}{-}: P \times A
\rightarrow R$ a bilinear form such that $\alpha_R$ is injective, then $P$ is
mock-projective.
\end{lemma}

\begin{proof}
Define
\begin{align*}
\underline{f} : S &\longrightarrow R^n \\ s &\longmapsto (f_1(s), \dots,
f_n(s)),
\end{align*}
and consider the $R$--submodule $M \subseteq R^n$ generated by the elements
$\underline{f}(s), s \in S$. Since $R$ is noetherian, $M$ is finitely generated,
and there are $s_1, \dots, s_m \in S$ for which $M = \sum R
\underline{f}(s_l)$. For any $s \in S$ we have the set
\begin{equation*}
X(s) = \{ (r_1, \dots, r_m) \in R^m ~|~ \underline{f}(s) = \sum r_l
\underline{f}(s_l) \} \neq \varnothing .
\end{equation*}
For each $s \in S$, choose $r(s) \in X(s)$. This gives a map $r : S \rightarrow
R^m$. Now. it is clear that there are maps $g_1, \dots, g_m \in R^S$ such that
$r = (g_1, \dots, g_m)$. Finally, $\underline{f}(s) = \sum g_l(s)
\underline{f}(s_l)$. But this is an equality in $R^n$, whence, for each $i = 1,
\dots, n$ we obtain $f_i(s) = \sum g_l(s) f_i(s_l)$, and hence $f_i =
\sum f_i(s_l) g_l$.
\end{proof}

\begin{remark}
Assume $\alpha_R$ to be injective. Let $M$ be an $R$--module. It
follows directly from Lemma \ref{movemaps} and Proposition
\ref{prop-Q} that $P \subseteq R^A$ is $M$--pure if and only if
$\alpha_M$ is injective.
\end{remark}

\begin{proposition}\label{dualQ}
Let $A$ be an $R$--algebra and assume $A^\circ$ is pure in $R^A$. Then
\begin{enumerate}
\item $A^\circ$ is an $R$--coalgebra. If in addition $A$ is a bialgebra (resp. Hopf algebra) then $A^\circ$ is a bialgebra (resp. Hopf algebra).
\item  Let $B$ be an $R$--algebra such that $B^\circ$ is pure in $R^B$. For every morphism of $R$--algebras $\varphi : A \rightarrow B$ we have
\begin{equation*}
\varphi^* (B^{\circ}) \subseteq A^{\circ},
\end{equation*}
moreover $\varphi^{\circ} := \varphi_{|_{B^{\circ}}}^*$ is an $R$--coalgebra
morphism.
\item Let $C \subseteq A^\circ$ be a subcoalgebra. Consider the bilinear form
\begin{equation*}
\begin{split}
\bil{-}{-} : C \times A &\longrightarrow R \\ (f,a) &\longmapsto \bil{f}{a} =
f(a).
\end{split}
\end{equation*}
Then $(C,A,\bil{-}{-})$ is a rational pairing.
\item If $(C,A,\bil{-}{-})$ is a rational pairing then $C$ is an
  $R$--subcoalgebra of $A^\circ$.
\end{enumerate}
\end{proposition}

\begin{proof}
(1) By \cite[Theorem 2.8]{Abuhlail/Gomez/Wisbauer:2000} $A^\circ$ is a
coalgebra. 
The rest of the first statement is similar to the argument in
\cite[9.1.3]{Montgomery:1993} due to the fact that over noetherian rings,
submodules of finitely generated modules are finitely generated.

(2) Let $f \in B^\circ$ and assume $I \subseteq B$ to be a cofinite left ideal
contained in $\ker f$. Since $R$ is noetherian, it is easy to check that
$\varphi^{-1}(I) \subseteq A$ is a left cofinite ideal contained in $\ker(f
\circ \varphi) = \ker(\varphi^*(f))$. A diagram chase shows that
$\varphi^\circ$ 
is an $R$--coalgebra map and the second statement is proved.

(3) Equation \eqref{eq:bilphi} is clearly satisfied. By
\cite[3.3]{Wisbauer:1997} $C$ is pure in $A^\circ$, so Proposition
\ref{prop-Q} and Lemma \ref{movemaps} give injectivity of $\alpha_M$ for each
$R$--module $M$.

(4) Since $\alpha_R$ is injective $C$ can be viewed as an $R$--submodule of
$A^*$. Let us see that $C \subseteq A^\circ$. If $c \in C$ and
$\Delta(c) = \sum c_1 \tensor c_2$, then $a c = \sum c_1 \bil{c_2}{a}$ by
Eq. \eqref{eq:bilphi} and \eqref{eq:dual-co}. It follows that $A c$ is
finitely generated by $\{ c_1 \}_{(c)}$ as $R$--module for every $c \in
C$. Moreover, Eq. \eqref{eq:bilphi} easily implies that the
comultiplication and the counit on $C$ are induced from $A^\circ$.
\end{proof}

\begin{corollary}\label{acs-pair}
Let $A$ be an $R$--algebra and assume $A^\circ$ is pure in $R^A$. Consider the
bilinear form
\begin{equation*}
\begin{split}
[-,-] : A^\circ \times {A^\circ}^* &\longrightarrow R \\ (f,\lambda)
&\longmapsto [f,\lambda] = \lambda(f)
\end{split}
\end{equation*}
for all $f \in A^\circ$ and $\lambda \in {A^\circ}^*$. Then
$(A^\circ,{A^\circ}^*,[-,-])$ is a rational pairing.
\end{corollary}

\begin{proof}
Assume $\sum [f_{i},\lambda] m_{i} = 0$ for all $\lambda
\in {A^\circ}^*$. Let $a \in A$ and consider
\begin{equation*}
\begin{split}
\bil{-}{a} : A^\circ &\longrightarrow R \\ f &\longmapsto \bil{f}{a} = f(a)
\end{split}
\end{equation*}
for all $f \in A^\circ$. Then $\sum [f_{i},\bil{-}{a}] m_{i} = 0$ and so
$\sum \bil{f_{i}}{a} m_{i} = 0$ for all $a \in A$, which implies $\sum
m_{i} \tensor f_{i} = 0$, since the pairing $(A^\circ,A,\bil{-}{-})$ is a
rational pairing by Proposition \ref{dualQ}. 
\end{proof}

As a consequence of Theorem \ref{iso}, Propositions
\ref{dualQ} and \ref{prop-q}, and Corollary \ref{acs-pair} we have:

\begin{theorem}\label{equal}
Let $A$ be an $R$--algebra such that $A^\circ$ is pure in $R^A$. Let $\varphi :
A \rightarrow A^{\circ *}$ be the canonical morphism and let $\varphi_* :
\lMod{A^{\circ *}} \rightarrow \lMod{A}$ be the restriction of scalars
functor. Then
\begin{enumerate}[(1)]
\item The functors $(-)^{A^{\circ}} : \Rat (\lMod{A}) \rightarrow
\rCom{A^{\circ}}$ and $(-)^{A^{\circ}} : \Rat (\lMod{A^{\circ *}}) \rightarrow
\rCom{A^{\circ}}$ are isomorphisms of categories.
\item $\Rat(\lMod{A}) = \mathbf{\sigma}[{}_{A}A^\circ]$ and $\Rat(\lMod{A^{\circ *}}) = \mathbf{\sigma}[{}_{A^{\circ
*}}A^\circ]$
\item The following diagram of functors is commutative
\begin{equation*}
\xymatrix{ \lMod{A^{\circ *}} \ar[d]_{\Rat^{T'}} \ar[rr]^{\varphi_*} & &
\lMod{A} \ar[d]^{\Rat^T} \\ \Rat^{T'}(\lMod{A^{\circ *}})
\ar[rd]_{(-)^{A^\circ}}^{\simeq} & & \Rat^T(\lMod{A})
\ar[ld]^{(-)^{A^\circ}}_{\simeq} \\ & {\rCom{A^\circ}} & }
\end{equation*}
where $T=(A^\circ,A,\bil{-}{-})$ and $T' = (A^{\circ},A^{\circ *},[-,-])$ are
the canonical pairings.
\end{enumerate}
\end{theorem}

\section{An application: Blattner--Montgomery duality.}

We are going to prove a Blattner--Montgomery like theorem (see
\cite[Theorem 9.4.9]{Montgomery:1993}) when \emph{$R$ is any commutative 
noetherian ring} and $(H,m,u,\Delta,\epsilon,S)$ is a Hopf algebra
such that $H^\circ$ is pure in $R^H$ (this condition holds if $H$
is $R$--projective and $H^{\circ}$ is pure in $H^*$). When $R$ is a
Dedekind domain, we obtain as a corollary the version given in
\cite{Cao-Yu/Nichols:1990}.

We are going to recall some definitions and notations. A left
$H$--module algebra is an $R$-algebra $(A,m_A,u_A)$ such that $A$
is a left $H$--module and $m_A,u_A$ are $H$--module maps. This
means in terms of Sweedler's notation:
\begin{equation*}
h(ab) = \sum (h_1 a)(h_2 b) \quad \text{and} \quad h 1_A = \epsilon(h)
1_A
\end{equation*}
Analogously $A$ is a right $H$--comodule algebra if $A$ is a right $H$--comodule
(via $\rho : A \rightarrow A \tensor H$) and $m_A,u_A$ are $H$--comodule maps,
i.e.,
\begin{equation*}
\rho(ab) = \sum a_0 b_0 \tensor a_1 b_1 \quad \text{and} \quad \rho(1_A) = 1_A
\tensor 1_H.
\end{equation*}

Let $A$ be a left $H$--module algebra where the $H$--module action is denoted by
$w_A : H \tensor A \rightarrow A$. The following composition of maps
\begin{equation*}
\xymatrix{ (A \tensor H) \tensor (A \tensor H) \ar@{-->}[dd]_-{m_{A \# H}}
\ar[rr]^-{id \tensor \Delta \tensor id^{\tensor 2}} && A \tensor H \tensor H
\tensor A \tensor H \ar[d]^-{\tau_{(34)}} \\ && A \tensor H \tensor A \tensor H
\tensor H \ar[d]^-{id \tensor w_A \tensor id^{\tensor 2}} \\ A \tensor H && A
\tensor A \tensor H \tensor H \ar[ll]^-{m_A \tensor m} }
\end{equation*}
provides a structure of an associative $R$--algebra on $A \tensor H$. This algebra
is called the smash product of $A$ and $H$, and it is denoted by $A \# H$. In
Sweedler's notation, the multiplication can be viewed as follows:
\begin{equation*}
(a \# h)(b \# k) = \sum a (h_1 b) \# h_2 k,
\end{equation*}
where $a \# h = a \tensor h$. Since $H^\circ$ is pure in $R^H$, by Proposition
\ref{dualQ} we have
$(H^\circ,\Delta^\circ,\epsilon^\circ,m^\circ,u^\circ,S^\circ)$ is a Hopf
algebra. The left (and right) action of $H$ on
$H^\circ$ described in \eqref{eq:dual-co} makes $H^\circ$ a left (and  right)
$H$--module algebra (see \cite[Example 4.1.10]{Montgomery:1993}). In order to
make the notation consistent with the literature we denote the left (resp.
right) action of $H$ on $H^\circ$ by $\rightharpoonup$ (resp. $\leftharpoonup$).
Let $U$ be a Hopf subalgebra of $H^\circ$ (by definition $U \subseteq H^\circ$
should be pure, see \cite[3.3]{Wisbauer:1997}). Then $U$ is also a left $H$--module algebra. The action can be
described as
\begin{equation*}
\xymatrix{ H \tensor U \ar[r]^-{id \tensor m^\circ} & H \tensor U \tensor U
\ar[r]^-{\tau_{(123)}} & U \tensor H \tensor U \ar[rr]^-{\bil{-}{-} \tensor
  id} && U, }
\end{equation*}
and in Sweedler's notation,
\begin{equation*}
h \rightharpoonup f = \sum f_1 \bil{f_2}{h}
\end{equation*}
which allows the construction of $U \# H$.

Analogously $H$ is a left (resp. right) $U$--module algebra via
\begin{gather*}
\xymatrix{ U \tensor H \ar[r]^-{id \tensor \Delta} & U \tensor H \tensor H
\ar[r]^-{\tau_{(23)}} & U \tensor H \tensor H \ar[rr]^-{\bil{-}{-} \tensor id}
&& H, } \\
\left( \text{resp.} \xymatrix{ H \tensor U \ar[r]^-{\Delta \tensor id} &
H \tensor H \tensor U \ar[r]^-{\tau_{(123)}} & U \tensor H \tensor H
\ar[rr]^-{\bil{-}{-} \tensor id} && H, } \right)
\end{gather*}
This action is denoted by $\rightharpoonup$ (resp. $\leftharpoonup$), and
Sweedler's notation means
\begin{equation*}
f \rightharpoonup h = \sum h_1 \bil{f}{h_2} \quad \left( \text{resp.} \quad h
\leftharpoonup f = \sum \bil{f}{h_1} h_2 \right)
\end{equation*}
and we can construct $H \# U$.

These actions and constructions are analogous to the ones over a field. See
\cite[1.6.5, 1.6.6, 4.1.10]{Montgomery:1993} for details. Following
\cite[Definition 9.4.1]{Montgomery:1993} we have the following maps:
\begin{align*}
\lambda: H \# U &\longrightarrow \operatorname{End}_R(H) \\ h \# f &\longmapsto
\big[ k \mapsto h (f \rightharpoonup k) \big] \\ \rho : U \# H &\longrightarrow
\operatorname{End}_R(H) \\ f \# h &\longmapsto \big[ k \mapsto (k \leftharpoonup
f) h \big]
\end{align*}

\begin{lemma}\label{9.4.2}
$\lambda$ is an algebra morphism and $\rho$ is an anti-algebra morphism. If also
$H$ has bijective antipode, then $\lambda$ and $\rho$ are injective.
\end{lemma}

\begin{proof}
Following \cite[Lemma 9.4.2]{Montgomery:1993}, we consider $\lambda$, as the
argument for $\rho$ is similar. Straightforward computations show that $\lambda$
is an algebra morphism. To see the injectivity we define $\lambda': H \# U
\rightarrow \operatorname{End}_R(H)$ and $\psi: \operatorname{End}_R(H)
\rightarrow \operatorname{End}_R(H)$ as follows:
\begin{gather*}
\lambda'(h \# f)(k) = \bil{f}{k}h \\ \psi(\sigma) = (\sigma \tensor
\overline{S}) \circ \tau \circ \Delta
\end{gather*}
where $\overline{S}$ is the composition inverse of $S$. We can see that
$\lambda' = \psi \circ \lambda$ as in \cite[Lemma 9.4.2]{Montgomery:1993}.
Moreover, $(U, H, \bil{-}{-})$ is a rational pairing by Proposition \ref{dualQ},
so $\lambda'$ is injective.
\end{proof}

We say that $U$ satisfies the RL--condition with respect to $H$ if $\rho(U \# 1)
\subseteq \lambda (H \# U)$.

Let $(A,\rho_A)$ be a right $U$--comodule algebra. Then $A$ is a left
$H$--module algebra with action
\begin{equation}\label{eq:action}
\xymatrix{ H \tensor A \ar[r]^-{id \tensor \rho_A} & H \tensor A \tensor U
\ar[r]^-{\tau_{(132)}} & A \tensor U \tensor H \ar[rr]^-{id \tensor
  \bil{-}{-}} && A },
\end{equation}
or in Sweedler's notation,
\begin{equation*}
h a = \sum a_0 \bil{a_1}{h}.
\end{equation*}

\begin{theorem}\label{BM-main}
Let $H$ be a Hopf algebra such that $H^\circ$ is pure in $R^H$, and let $U$ be
a Hopf subalgebra of $H^\circ$. Assume that both $H$ and $U$ have bijective
antipodes and $U$ satisfies the RL--condition with respect to $H$. Let $A$ be
a right $U$--comodule algebra. Let $U$ act on $A \# H$ by acting trivially on
$A$ and via $\rightharpoonup$ on $H$. Then
\begin{equation*}
(A \# H) \# U \simeq A \tensor (H \# U)
\end{equation*}
\end{theorem}

\begin{proof}
The computations in \cite[Theorem 9.4.9 and Lemma 9.4.10]{Montgomery:1993}
remain valid here once we have proved Lemma \ref{9.4.2}.
\end{proof}

\begin{remark}
Let $R$ be a Dedekind domain and assume that $A$ is an $U$--locally
finite left $H$--module algebra and that $U$ is $R$--closed in
$H^{\circ}$ in the sense of \cite{Cao-Yu/Nichols:1990}. By
\cite[Lemma 4]{Cao-Yu/Nichols:1990}, $A$ is a rational left
$H$--module which implies, by Theorem \ref{equal}, that $A$ is a
right $U$--comodule algebra. Therefore, \cite[Theorem
5]{Cao-Yu/Nichols:1990} follows as a corollary of Theorem
\ref{BM-main}.
\end{remark}

\begin{remark}
If $H$ is cocommutative then $U$ satisfies the RL--condition (see \cite[9.4.7
Example]{Montgomery:1993}), so examples in subsections \ref{Rx1...xn} and
\ref{laurent} and Example \ref{graduado} satisfy the RL--condition. So let $G$
be a group such that $R[G]^\circ$ is pure in $R^{R[G]}$ (if $G$ is either
finite or $R$ is hereditary, this condition is satisfied), and
let $A$ be an $R$--algebra such that $G$ acts as automorphisms on $A$. Then we have
\begin{displaymath}
(A \# R[G]) \# R[G]^\circ \cong A \tensor (R[G] \# R[G]^\circ)
\end{displaymath}
\end{remark}

\section{Examples}

\emph{In this section $R$ is assumed to be noetherian}. We are going to consider a
class of $R$--algebras for which $A^\circ$ is pure in $R^A$ (and hence $A^\circ$
has a structure of an $R$--coalgebra). For every $R$--algebra $A$ let
$\mathcal{L}_{\mathrm{cof}}$ be the linear topology on $A$ whose basic
neighborhoods of $0$ are the $R$--cofinite left ideals, i.e.,
\begin{equation*}
\mathcal{L}_{\mathrm{cof}} = \{I \leq {}_A A ~|~ A/I \text{ is finitely
generated as an $R$--module}\}
\end{equation*}

\subsection{The category $\PellAlg{}$}

\begin{definition}[{Property \Pleft{}}]
An $R$--algebra $A$ has property \Pleft{} in case the set
\begin{displaymath}
\mathcal{P}_{\mathrm{cof}} = \{I \leq {}_A A ~|~ A/I \text{ is finitely
  generated and projective as an $R$--module}\}
\end{displaymath}
is a basis for $\mathcal{L}_{\mathrm{cof}}$, i.e. for every left cofinite ideal
$I$ of $A$, there exists a left ideal $I_{0} \subseteq I$, with $A/I_{0}$
finitely generated and projective as an $R$--module.

We denote by $\PellAlg$ the full subcategory of $\mathrm{Alg}_{R}$ whose
objects are all $R$--algebras which have property \Pleft{}.
\end{definition}

\begin{proposition}\label{pro-co}
If $A \in \PellAlg$, then $(A^\circ,A,\bil{-}{-})$ is a rational system.
\end{proposition}

\begin{proof}
Let $M \in \lMod{R}$ and let $\sum m_i \tensor f_i \in M \tensor A^\circ$.
Assume $\sum m_i \bil{f_i}{a} = 0$ for every $a \in A$. Notice that for each
$i$, $f_i \in (A/I_i)^*$ for some cofinite ideal $I_i$ of $A$. Put $J =
\bigcap_i I_i$. Then $J$ is cofinite. Since $A \in \PellAlg$ there exists some
ideal $J_0 \subseteq J$ such that $A/J_0$ is finitely generated and projective
as $R$--module (and so $(A/J_0)^{**} \simeq A/J_0$). Let $\{ a_\lambda + J_0,
\phi_\lambda \}_\Lambda$ be a finite dual basis for $(A/J_0)^*$. Since $f_i \in
(A/J_0)^*$ for all $i$, we get
\begin{equation*}
\begin{split}
\sum m_i \tensor f_i &= \sum m_i \tensor \sum \bil{f_i}{a_\lambda + J_0}
\phi_\lambda \\
&= \sum \big( \sum \bil{f_i}{a_\lambda + J_0} m_i \big) \tensor \phi_\lambda
\\ 
&= \sum 0 \tensor \phi_\lambda = 0 \quad
\text{(notice that $f_i(J_0) = 0$)}.
\end{split}
\end{equation*}
Hence $(A^\circ,A,\bil{-}{-})$ is a rational system by Proposition \ref{prop-q}.
\end{proof}

\begin{corollary}\label{Pl=>coal}
If $A \in \PellAlg$ then
\begin{enumerate}
\item $A^\circ$ is an $R$--coalgebra. If in addition $A$ is a bialgebra (resp. Hopf algebra) then $A^\circ$ is a bialgebra (resp. Hopf algebra).
\item $(A^\circ,A,\bil{-}{-})$ is a rational pairing.
\end{enumerate}
\end{corollary}

\begin{proof}
It follows directly from Propositions \ref{prop-Q}, \ref{dualQ} and
\ref{pro-co}.
\end{proof}

\begin{remark}
By \eqref{eq:dual-co2}, the proof of Proposition \ref{pro-co} remains true if
we replace left ideals in property \Pleft{} by right or two sided ones. So we
can speak of property $\mathbf{P_r}$ or property $\mathbf{P}$.
\end{remark}

\begin{remark}
If $A \in \PellAlg$, then
\begin{equation*}
\begin{split}
A^{\circ *} &= \ModHom{R}(A^\circ,R) = \ModHom{R} (\varinjlim_{I \in
  \mathcal{L}_{\mathrm{cof}}}(A/I)^*,R) \\ 
&\simeq \ModHom{R} (\varinjlim_{I \in \mathcal{P}_{\mathrm{cof}}}(A/I)^*,R) \\ 
&\simeq \varprojlim_{I\in \mathcal{P}_{\mathrm{cof}}}(A/I)^{**} \simeq
\varprojlim_{I\in \mathcal{P}_{\mathrm{cof}}} A/I \\ 
&\simeq \varprojlim_{I\in \mathcal{L}_{\mathrm{cof}}} A/I = \hat{A},
\end{split}
\end{equation*}
which means that $A^{\circ *}\simeq \hat{A}$, the completion of $A$ with
respect to the cofinite topology.
\end{remark}

\begin{proposition}\label{shiftPleft}
Let $A$ be in $\PellAlg$ and let $B$ be an $R$-algebra
extension of $A$ such that $B$ is finitely generated and projective in
$\rMod{A}$. Then $B$ belongs to $\PellAlg$.
\end{proposition}

\begin{proof}
Let $J \leq B$ be a cofinite left ideal. Then $J \cap A \leq A$ is a cofinite
left ideal because $R$ is noetherian. Since $A$ belongs to $\PellAlg$, there exists $I_0
\subseteq J \cap A$ such that $A/I_0$ is finitely generated and projective in
$\lMod{R}$. By the natural isomorphism
\begin{equation*}
\ModHom{R} \left( B \tensor[A] \tfrac{A}{I_0}, - \right) \cong \ModHom{A} \left(
B , \ModHom{R} \left( \tfrac{A}{I_0} , - \right) \right)
\end{equation*}
$B \tensor[A] \frac{A}{I_0}$ is finitely generated and projective in $\lMod{R}$.
Since $BI_0 \subseteq J$ and  $B \tensor[A] \frac{A}{I_0} \cong \frac{B}{BI_0}$
we get $B$ is in $\PellAlg$.
\end{proof}

\begin{example}\label{graduado}
Let $G$ be a group. An $R$--algebra is called $G$--graded if for every $\sigma
\in G$ there exists an $R$--submodule $A_\sigma \subseteq A$ such that $A =
\bigoplus_{\sigma \in G} A_\sigma$ and $A_\sigma A_\tau \subseteq
A_{\sigma\tau}$. If in addition $A_\sigma A_\tau = A_{\sigma\tau}$, $A$ is
called strongly graded. Let $G$ be finite with neutral element $e$ and let $A$
be a strongly $G$--graded $R$--algebra. By \cite[I.3.3
Corollary]{Nastasescu/Oystaeyen:1982} it is clear that $A$ is finitely generated
and projective as right $A_e$--module, so if $A_e$ is in $\PellAlg$ then $A$ also belongs to $\PellAlg$. In particular, if $A$ is in $\PellAlg$ and $G$ is a finite group, then
a crossed product $A * G$ also belongs to $\PellAlg$. Crossed products cover
the following cases: if $A \in \PellAlg$ then $A[G], A^t[G], AG \in \PellAlg$
where $A[G]$ is the group algebra, $A^t[G]$ is the twisted group algebra and
$AG$ is the skew group algebra. See \cite{Passman:1989} for an introduction on
crossed products.
\end{example}

Our aim is the proof of Theorem \ref{main}, which was shown in \cite[Lemma
6.0.1.]{Sweedler:1969} for algebras over fields and in
\cite{Cao-Yu/Nichols:1990} for algebras over Dedekind domains. However we need
some technical statements.

\begin{lemma}\label{tensor}
Let $M$ and $N$ be two $R$--modules and consider submodules $M' \subseteq M$
and $N' \subseteq N$. Assume $M'$ to be $N$--pure and $N'$ to be $M$--pure (this
is in particular valid if $M$ and $N$ are flat in $\lMod{R}$). Then
\begin{equation*}
M/M' \tensor N/N' \simeq (M \tensor N)/(M' \tensor N + M \tensor N').
\end{equation*}
\end{lemma}

\begin{proof}
By purity $M' \tensor N$ and $M \tensor N'$ are $R$--submodules of $M \tensor
N$. Since the diagram
\begin{equation*}
\xymatrix{ M \tensor N \ar@{->>}[r] \ar@{->>}[d] & M/M' \tensor N \ar@{->>}[d]
\\ M \tensor N/N' \ar@{->>}[r] & M/M' \tensor N/N' }
\end{equation*}
is a pushout diagram, the result follows.
\end{proof}

\begin{proposition}\label{AoB}
Let $A,B$ be algebras in $\PellAlg$.
\begin{enumerate}
\item If $K \leq A \tensor B$ is a cofinite left ideal then there exist $I_0 \leq A$ and $J_0 \leq B$ such that $A/I_0$ and $B/J_0$ are finitely generated and projective in $\lMod{R}$, and so that $I_0 \tensor B + A \tensor J_0 \subseteq K$.
\item The $R$--algebra $A \tensor B$ belongs to $\PellAlg$.
\end{enumerate}
\end{proposition}

\begin{proof}
Consider the canonical maps
\begin{equation*}
\begin{split}
\alpha:A &\longrightarrow A \tensor B \\ a &\longmapsto a \tensor 1
\end{split}
\end{equation*}
and
\begin{equation*}
\begin{split}
\beta:B &\longrightarrow A \tensor B \\ b &\longmapsto 1 \tensor b
\end{split}
\end{equation*}
Put $I=\alpha^{-1}(K)$ and $J=\beta^{-1}(K)$. Since $R$ is noetherian, $I$ and
$J$ are cofinite left ideals of $A$ and $B$, respectively. Since $A, B \in
\PellAlg$ there exist $I_0 \subseteq I$ and $J_0 \subseteq J$ such that $A/I_0$
and $B/J_0$ are finitely generated and projective in $\lMod{R}$. Let $K_0 = I_0
\tensor B + A \tensor J_0$. Since $I_0 \leq A$ and $J_0 \leq B$ are pure
submodules we have $K_0 \subseteq K$ as desired.

By Lemma \ref{tensor}
\begin{equation*}
\frac{A \tensor B}{K_0} \simeq \frac{A}{I_0} \tensor \frac{B}{J_0},
\end{equation*}
hence $(A \tensor B)/K_0$ is finitely generated and projective and $A \tensor B$ is in $\PellAlg$.
\end{proof}

\begin{theorem}\label{main}
Let $A,B$ be in $\PellAlg$. Then there is a canonical isomorphism $A^{\circ} \tensor B^{\circ} \simeq (A \tensor B)^{\circ}$.
\end{theorem}

\begin{proof}
Since $A,B$ are in $\PellAlg$, $A^\circ$ is pure in $R^A$ and $B^\circ$ is pure in
$R^B$ by Propositions \ref{pro-co} and \ref{prop-Q}. So $A^\circ
\tensor B^\circ \subseteq R^A \tensor R^B$. Let $\pi$ be the morphism:
\begin{equation*}
\begin{split}
\pi : R^A \tensor R^B &\longmapsto R^{A \times B} \\ f \tensor g &\longmapsto
\left[ (a,b) \mapsto f(a)g(b) \right]
\end{split}
\end{equation*}
By \cite[Proposition 1.2]{Abuhlail/Gomez/Wisbauer:2000} this map is injective, so
the statement will be clear once we have seen $\pi(A^\circ \tensor B^\circ) = (A
\tensor B)^\circ$. So let $f \tensor g \in A^\circ \tensor B^\circ$ and let $I
\subseteq A$ and $J \subseteq B$ be left ideals contained in $\ker f$ and $\ker
g$ respectively and such that $A/I$, $B/J$ are finitely generated and projective
(they exist because $A,B$ belong to $\PellAlg$). Since $I \subseteq A$ and $J \subseteq
B$ are pure, by Lemma \ref{tensor} $I \tensor B + A \tensor J \subseteq A
\tensor B$ is a cofinite left ideal, which is contained in $\ker(\pi(f \tensor
g))$. As $\pi(f \tensor g)$ is bilinear it is clear that $\pi(A^\circ \tensor
B^\circ) \subseteq (A \tensor B)^\circ$.

Let $h \in (A \tensor B)^\circ$, and assume $K \subseteq A \tensor B$ to be a
cofinite left ideal contained in $\ker h$. By Proposition \ref{AoB} there exist left
ideals $I_0 \leq A$ and $J_0 \leq B$ such that $A/I_0$ and $B/J_0$ are finitely
generated and projective in $\lMod{R}$ and so that $I_0 \tensor B + A \tensor
J_0 \subseteq K$. By Lemma \ref{tensor} there is an epimorphism
\begin{equation*}
\frac{A}{I_0} \tensor \frac{B}{J_0} \rightarrow \frac{A \tensor B}{K}
\rightarrow 0
\end{equation*}
which induces a monomorphism
\begin{equation*}
0 \rightarrow \left( \frac{A \tensor B}{K} \right)^* \rightarrow \left(
\frac{A}{I_0} \tensor \frac{B}{J_0} \right)^* \simeq \left( \frac{A}{I_0}
\right)^* \tensor \left( \frac{B}{J_0} \right)^* \subseteq A^\circ \tensor
B^\circ
\end{equation*}
So there exist elements $f_1, \dots, f_n \in (A/I_0)^* \subseteq A^\circ$ and
$g_1, \dots, g_n \in (B/J_0)^* \subseteq B^\circ$ such that $\pi(\sum f_i \tensor g_i) = h$. This completes the proof.
\end{proof}

We finish with some examples.

\subsection{The $R$--bialgebra $R[x_{1},...,x_{n}]^\circ$}\label{Rx1...xn}\label{rx-coal}\label{rxk-bi}

By \cite[Proposition 3.1]{Abuhlail/Gomez/Wisbauer:2000}, every cofinite ideal
$I \leq R[x]$ contains a monic polynomial $f(x)$). Put $I_{0}=(f(x)) \subseteq
I$. Then $R[x]/I_{0}$ is finitely generated and projective (in fact
free). Hence $R[x]$ is in $\PellAlg$ and so $R[x]^\circ$ is an $R$--coalgebra
by Corollary \ref{Pl=>coal}. Moreover, $R[x_{1},\dots ,x_{n}]$ belongs to
$\PellAlg$ by Proposition \ref{AoB}. There are two canonical bialgebra
structures on $R[x_1, \dots, x_n]$. The first one comes from the semigroup
algebra structure of $R[x_1, \dots, x_n]$ (i.e. every $x_i$ is a group-like
element), and the second one appears when we see $R[x_1, \dots, x_n]$ as the
enveloping algebra of an abelian Lie algebra (i.e. every $x_i$ is a primitive
element). The latter one 
is a Hopf algebra structure. By Corollary \ref{Pl=>coal}, $R[x_1, \dots,
x_n]^\circ$ is a bialgebra (resp. Hopf algebra).

It follows from Proposition \ref{AoB} that if $A$ belongs to $\PellAlg$ then $A[x_1, \dots, x_n]$ is in $\PellAlg$.

\subsection{The Hopf $R$--algebra of Laurent polynomials}\label{laurent}

\begin{definition}
A monic polynomial $q(x)\in R[x]$ is called reversible if $q(0)$ is a unit in
$R$. An ideal $I \subseteq R[x,x^{-1}]$ is called reversible if it contains a
reversible polynomial $q(x)$.
\end{definition}

\begin{lemma}\label{q-q}
Let $q(x) \in R[x]$ be a reversible polynomial. Then
\begin{equation*}
R[x]/(q(x)) \simeq R[x,x^{-1}]/(q(x)).
\end{equation*}
\end{lemma}

\begin{proof}
Let $q(x)=x^{n}+a_{n-1}x^{n-1}+\cdots +a_{1}x+a_{0}$ be a reversible
polynomial (i.e. $a_{0}$ in a unit in $R$). Notice
\begin{equation*}
R[x,x^{-1}]/(q(x)) \simeq R[x,y]/(xy-1,q(x)).
\end{equation*}
Put $I=(xy-1,q(x))$ and consider the $R$--linear map
\begin{equation*}
\begin{split}
\Psi : R[x] &\longrightarrow R[x,y]/I \\ x &\longmapsto x + I.
\end{split}
\end{equation*}
Clearly $\Psi$ is an $R$--algebra homomorphism and $\ker(\Psi)=I \cap R[x]$.
Moreover, $I \cap R[x] = R[x]q(x)$. Clearly $\Psi$ is surjective if and only if
$y+I \in \im(\Psi)$. Notice
\begin{equation*}
\begin{split}
yq(x)-x^{n-1}(xy-1) &= a_{0}y+a_{1}yx+\cdots+a_{n-1}yx^{n-1}+x^{n-1} \\ &=
a_{0}y+a_{1}+a_{2}x+\cdots +a_{n-1}x^{n-2}+x^{n-1} \bmod(I)
\end{split}
\end{equation*}
So
\begin{equation*}
y = -a_{0}^{-1}\,[x^{n-1}+a_{n-1}x^{n-2}+\cdots +a_{1}] \bmod(I)
\end{equation*}
Hence $y \in \im(\Psi)$ and we conclude that $\Psi$ is surjective.
\end{proof}

\begin{proposition}\label{m-ideal}
\begin{enumerate}
\item Let $I\subseteq R[x,x^{-1}]$ be a reversible ideal. Then $R[x,x^{-1}]/I$ is finitely generated as an $R$--module.
\item Let $R$ be noetherian. Assume $R[x,x^{-1}]/I$ to be finitely generated as an $R$--module. Then $I$ is reversible ideal.
\end{enumerate}
\end{proposition}

\begin{proof}
(1)  Let $I \subseteq R[x,x^{-1}]$ be a reversible ideal. Then $I$ contains a
reversible polynomial $q(x)$. By Lemma \ref{q-q}, $R[x,x^{-1}]/(q(x)) \simeq
R[x]/(q(x))$ which implies, by \cite[Proposition
3.1]{Abuhlail/Gomez/Wisbauer:2000}), that $R[x,x^{-1}]/(q(x))$ is finitely
generated as an $R$--module. Therefore, $R[x,x^{-1}]/I$ is finitely generated as
an $R$--module.

(2) Since $R[x]/(R[x]\cap I)$ embeds in the finitely generated $R$--module
$R[x,x^{-1}]/I$, we get that $R[x]/(R[x] \cap I)$ is finitely generated as an
$R$--module. By \cite{Abuhlail/Gomez/Wisbauer:2000}, there exists a monic
polynomial $f_1(x) = a_0 + a_1 x + \dots + x^n \in I \cap R[x]$. We know
$R[x,x^{-1}]$ is a Hopf $R$--algebra with antipode
\begin{equation*}
\begin{split}
S:R[x,x^{-1}] &\longrightarrow R[x,x^{-1}] \\ x &\longmapsto x^{-1},
\end{split}
\end{equation*}
Since $S$ is bijective, $R[x,x^{-1}]/I \simeq R[x,x^{-1}]/S(I)$ as $R$--modules.
So there exists a monic $f_2(x) = b_0 + \dots + b_{m-1} x^{m-1} + x^m \in S(I)
\cap R[x]$. Hence we have that $q(x) = x^m(f_1(x) + S(f_2(x))) \in I$. An easy
computation shows that
\begin{equation*}
q(x) = 1 + b_{m-1} x + \dots + (b_0 + a_0) x^m + \dots a_{n-1} x^{n+m-1} +
x^{n+m}
\end{equation*}
and $I$ contains the reversible polynomial $q(x)$. By Lemma \ref{q-q}
\begin{equation*}
R[x,x^{-1}]/(q(x)) \simeq R[x]/(q(x))
\end{equation*}
and so is finitely generated and projective (in fact free) as an $R$--module.
\end{proof}

\begin{theorem}\label{haupt}
Let $R$ be noetherian and let $A$ be in $\PellAlg$. Then $A[x_1, x_1^{-1}, \dots, x_n, x_n^{-1}]$ belongs to $\PellAlg$. In particular $R[x_1,x_1^{-1}, \dots, x_n, x_n^{-1}]^\circ$ is a Hopf algebra.
\end{theorem}

\begin{proof}
By Proposition \ref{m-ideal} and Lemma \ref{q-q} it is easy to see
that $R[x,x^{-1}]$ is in $\PellAlg$, so the first statement follows
from Proposition \ref{AoB}. Since $R[x_1,x_1^{-1}, \dots, x_n,
x_n^{-1}]$ is a group algebra, the last assertion follows from
Corollary \ref{Pl=>coal}.
\end{proof}

\bibliographystyle{amsplain}

\end{document}